\documentstyle[12pt]{article}
 
\begin{document}
  
 \def \wt{{\rm wt}}
\def \tr{{\rm tr}}
\def \span{{\rm span}}
\def \Res{{\rm Res}}
\def \Res{{\rm QRes}}
\def \End{{\rm End}}
\def \E{{\rm End}}
\def \Ind {{\rm Ind}}
\def \Irr {{\rm Irr}}
\def \Aut{{\rm Aut}}
\def \Hom{{\rm Hom}}
\def \mod{{\rm mod}}
\def \ann{{\rm Ann}}
\def \<{\langle} 
\def \>{\rangle} 
\def \t{\tau }
\def \a{\alpha }
\def \e{\epsilon }
\def \l{\lambda }
\def \L{\Lambda }
\def \b{\beta }   
\def \om{\omega }
\def \o{\omega }
\def \c{\chi}
\def \ch{\chi}
\def \cg{\chi_g}
\def \ag{\alpha_g}
\def \ah{\alpha_h}
\def \ph{\psi_h}
\def \be{\begin{equation}}
\def \ee{\end{equation}}
\def \bex{\begin{example}\label}
\def \bea{\begin{eqnarray} }
\def \eea{\end{eqnarray} }
\def \eex{\end{example}}
\def \bl{\begin{lem}\label}
\def \el{\end{lem}}
\def \bt{\begin{thm}\label}
\def \et{\end{thm}}
\def \bp{\begin{prop}\label}
\def \ep{\end{prop}}
\def \br{\begin{rem}\label}
\def \er{\end{rem}}
\def \bc{\begin{coro}\label}
\def \ec{\end{coro}}
\def \bd{\begin{de}\label}
\def \ed{\end{de}}
\def \pf{{\bf Proof. }}
\def \voa{{vertex operator algebra}}

\newtheorem{theorem}{Theorem}[section]
\newtheorem{corollary}{Corrolary}[section]
\newtheorem{lemma}{Lemma}[section]
\newtheorem{conjecture}{Conjecture}[section]
\newtheorem{convention}{Convention}[section]
\newtheorem{remark}{Remark}[section]
\newtheorem{definition}{Definition}[section]
\newtheorem{proposition}{Proposition}[section]
\renewcommand{\baselinestretch}{1.4}
\newcommand{\binom}{ { { {\rm {wt}}  a} \choose i } }
\newcommand{\wta}{{\rm {wt} }  a } 
\newcommand{\R}{\mbox {Res} _{z} }
\newcommand{\la}{\langle}
\newcommand{\ra}{\rangle}
\newcommand{\wtb}{{\rm {wt} }  b }  
\newcommand{\g}{\bf g}
\newcommand{\hg}{\hat {\bf g} }
\newcommand{\h}{\bf h}
\newcommand{\hh}{\hat {\bf h} }
\newcommand{\nn}{\bf n}
\newcommand{\Z}{\bf Z}
\newcommand{\N}{\bf N}
\newcommand{\C}{\bf C}
\newcommand{\Q}{\bf Q}
\newcommand{\1}{\bf 1}
\newcommand{\NS}{\bf NS}
\newcommand{\cA}{\cal A}
\newcommand{\hn} {\hat {\bf n}}
\newcommand{\hf}{\frac{1}{2} }   
\newcommand{\non}{\nonumber}
\def \l {\lambda}
\baselineskip=12pt
\newcommand{\rw}{\rightarrow}
\newcommand{\n}{\:^{\circ}_{\circ}\:}
\makeatletter
\@addtoreset{equation}{section}
\def\theequation{\thesection.\arabic{equation}}
\makeatother
\makeatletter

\title {Representations of $N=2$ superconformal
vertex algebra }
 
\author{ Dra\v zen Adamovi\' c }
\date{}
\pagestyle{myheadings}
\markboth
{Dra\v{z}en Adamovi\'{c}}
{}                          
\maketitle

\section*{Introduction}

In the last few years $N=2$ superconformal algebra has attracted
very much interests. The most important step in the representation
theory of $N=2$ superconformal algebra was made in the series of papers
\cite{FST}, \cite{ST}, \cite{FSST}. It was proved that certain categories of modules of
$N=2$ superconformal algebra and affine Lie algebra $A_1 ^ {(1)}$ are
equivalent [FST]. The structure of the singular vectors in highest
weight representations and embedding diagrams for Verma modules is
also very much understood (cf. \cite{D}, \cite{ST}).

In this paper we will investigate the representation theory of $N=2$
superconformal algebra from the vertex algebra point of view.
On the irreducible highest weight module $L_c$ exists a natural
structure of vertex operator superalgebra (SVOA) (cf. Section 1).
We will classify all irreducible $L_c$--modules.
The problem of the classification of irreducible $L_c$--modules
was initiated by Eholzer and Gaberdiel in \cite{EG}. They proved that if
$L_c$ is rational SVOA, then $L_c$ has to be unitary representations,
i.e. $c = \frac{3 m}{m+2}$ for $m \in {\N}$.
They used the  theory  of Zhu's algebra $A(L_c)$, and calculated it in some
special cases. Unfortunately, the complicated structure of singular vectors makes
the explicit determination of Zhu's algebra $A(L_c)$ extremely difficulty,
and they didn't get the complete classification result.

Instead of explicit calculation of Zhu's algebra, for problem of the classification
of irreducible $L_c$--modules we use the equivalence of categories from
\cite{FST}.  
We interpret the Kazama-Suzuki and anti Kazama-Suzuki mapping from \cite{FST}
in the language of vertex algebras and get embedding results between
certain simple vertex operator (super)algebras (cf. Sections 4, 5).
In this way we can use the representation theory of VOA $L(m,0)$ associated
to admissible $A_1 ^ {(1)}$--representation from \cite{AM}. As a result
we get the complete classification of irreducible $L_c$--modules.
     
Let us explain the classification result in more details. 
For $m = \frac{t}{u}$ admissible set $c_m =\frac{3m}{m+2}$,
$N= 2 u + t -2$ and
$$
S^m = \{ n - k (m+2) \ \vert \ k,n \in {\Z}_+, n\le N, k\le u-1\}.
$$
Let      
$$
W^{ c_m} = \left \{ \left( \frac{j k - \frac{1}{4} }{m+2},
 \frac{j-k}{m+2} \right) \ \vert \ j,k \in {\N} _{ \frac{1}{2} },
  0< j,k, j+k \le N+1 \right \},
$$ 
and if $m \notin {\N}$ let
$$
D^{c_m} =
\{ (h,q) \in {\C} ^2 \ \vert \ q^2 + \frac{4h}{m+2} = 
\frac{r (r+2)}{ (m+2) ^2}, \ r\in S^m \setminus {\Z} \}.
$$
Note that the set $W ^{ c_m}$ is finite and the set $D^{c_m}$
is union of finitely many rational curves.

Let $L_{h,q,c}$ denotes the irreducible highest weight $N=2$--module 
with the highest weight $(h,q,c)$ (see Section 1.).

\begin{theorem}  \label{unitary.rat}
Let $m \in {\N}$. Then  the set
$$
\{ L_{h,q,c_m} \ \vert \ (h,q) \in W^{c_m} \}
$$
provides all irreducible $L_{ c_m}$--modules.
\end{theorem}
If $m \in {\N}$ then  
$ L_{h,q,c_m}$, $(h,q) \in W^{c_m}$ are all unitary representations with  
the central charge $c_m$. So, in that case the irreducible $L_{c_m}$
modules coincide with the unitary representations for $N=2$ superconformal
algebra.
                                                                            
\begin{theorem}
Let $m \in {\Q}$ is admissible and $ m \notin {\N}$.
Then the set
$$
\{ L_{h,q,c_m} \ \vert \ (h,q) \in W^{ c_m} \cup D^{c_m} \}
$$
provides all irreducible $L_{c_m}$--modules.
\end{theorem}
 SVOA $L_{c_m}$, $m \notin {\N}$, has
 uncountably many non-isomorphic irreducible
 representations which have  
 a nice description as a union of finite set $W^{ c_m}$
and the infinite set $D ^{c_m}$ which is 
described with finitely many rational curves.

\section{$N=2$ superconformal vertex algebra}

In this section we will show that on the vacuum representations
of $N=2$ superconformal algebra exists a natural structure of vertex
operator superalgebra. This result follows from the results on local
generating fields of SVOAs (cf. \cite {K}, \cite{Li}, \cite{P}).
 This structure was alredy studied in \cite{EG}.

$N=2$ superconformal algebra ${\cA}$
is the infinite-dimensional Lie super algebra  with basis
$ L_n, T_n, G_r ^{\pm}, C$, $n\in {\Z}$, $ r\in {\hf} + {\Z}$ 
and (anti)commutation relations
given by
\bea
&& [L_m,L_n] = (m-n) L_{m+n} + \frac{C} {12} (m ^ 3 -m) \delta_{m+n,0}
\nonumber \\
&& [L_m, G_r ^{\pm}] =( {\hf } m -r )    G_{m+r} ^{\pm}  \non \\
&&[L_m, T_n] = - n T_{n+m} \non \\
&& [T_m, T_n ] =  \frac {C}{3} m \delta _{m+n,0} \non \\
&&  [T_m, G_r ^{\pm} ] = \pm G_{m+r} ^{\pm} \non \\
&& \{ G_r ^{+}, G_r ^{- }  \} = 2 L_{r+s} + (r-s) T_{r+s} +
    \frac{C}{3} ( r  ^ 2 - \frac{1}{4} ) \delta_{r+s,0} \non \\
&& [L_m,C]= [T_n,C] = [G_r ^{\pm}, C] =0 \non \\
&& \{ G_r ^{+}, G_s ^{+} \}= \{ G_r ^{-}, G_s ^{-} \} = 0 \non
\eea
for all $m,n \in  {\Z}$, $r,s \in   {\hf} + {\Z}$.

We denote the Verma module generated from a highest weight vector
$\vert h,q,c \rangle$
with $L_0$ eigenvalue $h$, $T_0$ eigenvalue $q$ and central charge $c$
by $ M _{h,q,c}$. An element $v \in M _{h,q,c}$ is called singular
vector if
$$
L_n v = T_n v = G_r ^{\pm} v = 0,
\ \ n, r+ {\hf} \in {\N},
$$
and $v$ is an eigenvector of $L_0$ and $T_0$.
Let $J _{h,q,c}$ be the maximal $U({\cal A})$--submodule in  $ M _{h,q,c}$.
Then
$$L_{h,q,c} =\frac{M _{h,q,c}}{ J _{h,q,c} }$$
is the irreducible highest weight module.

Now we will consider the Verma module $M_{0,0,c}$.
One easily sees that  for every $c \in {\C}$ 
$$G_{ {-\hf} } ^{\pm} \vert 0,0,c \rangle $$
are the singular vectors in $M_{0,0,c}$.
Set
$$
V_c = \frac {M_{0,0,c} }
{ U( {\cA} ) G_{ {-\hf} } ^{+} \vert 0,0,c \rangle   +
 U( {\cA} ) G_{ {-\hf} } ^{-} \vert 0,0,c \rangle        } .
  $$
Then $V_c$ is a highest weight  ${\cA}$--module. Let ${\1}$ denote the
highest weight vector. Let $L_c=L_{0,0,c}$ be the corresponding simple module.
Define the following four vectors in $V_c$:
$$
\tau ^{\pm} = G_{-\frac {3}{2} } ^{\pm} {\1}, \ 
 j = T_{-1} {\1} , \ 
  \nu = L_{-2} {\1}, 
$$
and set
\bea
&& G^{+} (z) = Y(\tau ^{+} ,z)
= \sum _{n \in {\Z} } G_{n+ {\hf} } ^{+} z ^{-n-2}, \non \\
&& G^{-} (z) = Y(\tau ^{-} ,z)
= \sum _{n \in {\Z} } G_{n+ {\hf} } ^{-} z ^{-n-2}, \non \\
&& L (z) = Y(\nu,z)
= \sum _{n \in {\Z} } L_{n }  z ^{-n-2}, \non \\
&& T (z) = Y(j,z)
= \sum _{n \in {\Z} } T_{n }  z ^{-n-1}. {\label{polja}}
\eea
It is easy to see that the fields  $G^{+} (z)$,
$G^{-} (z)$, $L(z)$, $T(z)$ are mutually local and the theory of local
fields (cf. \cite{K}, \cite{Li}, \cite{P} ) implies the following result.

\begin{proposition}
There is a unique exstension of the fields (\ref{polja}) such that
$V_c$ becomes vertex operator superalgebra (SVOA).  Moreover, $L_c$ is
a simple SVOA.
\end{proposition}



 \begin{definition} Let $V$ be SVOA. We will say that $V$ is $N=2$
 SVOA if there exist vectors $\tau ^{\pm}, \nu, j \in V$  such that
 components of the fields $Y(\tau ^{\pm},z)$, $Y(\nu,z)$, $Y(j,z)$  span
  $N=2$ superconformal algebra.
 \end{definition}

 Previous definition has the following obvious  but important
 consequence.
 \begin{corollary} \label{n2con}
 Assume that $V$ is $N=2$ SVOA. Then we have the following:
 \item[(1)]
 $V$ is a $U(\cal A)$--module.
 \item[(2)]
 $V$ is a module for the SVOA $V_c$.
 \item[(3)]
  $U(\cal A) . {\1}$ is a subalgebra of $V$, isomorphic to
  $V_c$ or to certain quotient of $V_c$.
\end{corollary}
  
  In what follows will present
 one construction of $N=2$ SVOA.

 \begin{remark}
 In \cite{Z}, Zhu constructed an associative algebra $A(V)$ for an arbitrary
 VOA $V$ and established a one to one correspondence between irreducible
 representations of $V$ and irreducible representation of $A(V)$. V. Kac
 and W. Wang extended in \cite{KWn} the definition of Zhu's algebra to
 the case of SVOAs. In our case, one can show that the Zhu's algebra
 $A(V_c)$ is isomorphic to the polynomial algebra ${\C}[x,y]$, and Zhu's
 algebra $A(L_c)$ is a certain quotient of  ${\C}[x,y]$.
 This will imply that every irreducible $L_c$--module
 has to be the irreducible highest weight $U({\cal A})$--module. Eholzer and
 Gaberdiel  in \cite{EG} starting from physical motivated definition of
 Zhu's algebra showed that for non-generic $c$
 $A(L_c)= {\C}[x,y] / I$, where $I$ is an ideal
 in   ${\C}[x,y]$ generated by two polinomials $p_1(x,y)$, $p_2(x,y)$. They
 also calculated $p_1$, $p_2$ for some special cases (see Table 3.1 in
 \cite{EG}).
 \end{remark}

\section {Vertex operator algebras associated to affine Lie algebra
${ \hat sl_2}$       }

In this section we recall the classification of the irreducible modules
for VOAs associated to affine Lie algebra $A_1 ^ {(1)}$ obtained by
the author and A. Milas in \cite{AM}.

Let ${\g}$ be a finite-dimensional simple Lie algebra over ${\C}$. 
The affine Lie algebra
${\hg}$ associated with ${\g}$ is defined as
$$
{\g} \otimes {\C}[t,t^{-1}] \oplus {\C}c
$$
with the usual commutation relations.   
Let  ${\g} = {\nn}_- + {\h} + {\nn}_+$ and  ${\hg}={\hn}_- +{\hh} + {\hn}_+$
be the usual triangular decompositions for ${\g}$ and ${\hg}$ .
Let $P = {\bf C} [t] \otimes {\g} \oplus {\bf C}c $ be upper parabolic 
subalgebra.
Let $U$ be any ${\g}$--module. Considering $U$ as a $P$--module,
we defined the induced module ( so called generalized Verma module)
$
M(m,U) = U({\hg} ) \otimes _{U(P)} U $, where the central element $c$ acts
as multiplication with ${m} \in {\C}$.
Clearly ${\hg}$--module $M(m,U)$ contains a maximal submodule which
intersect $U$ trivialy. Let $L(m, U)$ be the corresponding quotient.
If $U$ is an irreducible ${\g}$--module, then $L(m, U)$ is
an irreducible ${\hg}$--module at the level ${m}$.
For $\lambda \in {\h} ^*$ with $V(\lambda) $ we denote the irreducible highest
weight ${\g}$--module. Set $M(m,\lambda)=M(m, V(\lambda))$. Let  
$L(m,\lambda)$ denotes its irreducible quotient. 

\begin{theorem} ([FZ]) Every $M(m,0)$ $m \ne -g$ (where 
$g$ denotes dual Coxeter number) has the structure of VOA. 
Let $U$ be any ${\g}$--module. Then
every $M(m, U)$ is a module for $M(m,0)$. In particulary $M(m, \lambda)$
is $M(m,0)$--module.  On the irreducible highest weeight ${\hg}$--module
$L(m,0)$ exists the structure of simple VOA.
\end{theorem}

   Let now ${\g} =  sl_2$, with generators $e,f,h$ and relations
   $[e,f]=2 h$, $[h,e]=2 e$, $[h,f] = -2f$.
   Let $\Lambda_0, \ \Lambda_1$ denote the fundamental weights for ${\hg}$,
   and let $\omega_1$ the fundamental weight for ${\g}$.
   Let $r \in {\C}$. Then
   $L(m, V(r \omega_1)) =  L(m,r \omega_1) $ is the irreducible highest weight
   ${\hg}$--module with the highest weight
   $(m -r)\Lambda_0 + r \Lambda_1$.

\begin{definition} A rational number $m=t/u$ is called admissible if
$u\in {\N}$, $t\in {\Z}$, $(t,u)=1$ and $2u+t-2 \ge 0$.
\end{definition}

Let $m=t/u \ \in {\Q}$ be admissible.  
We define the following set of weights at the level $m$: 
\bea    
P^{m} = \{ {\l}_{m,k,n}=
(m-n+k (m+2) ){\Lambda_0} + (n-k(m+2)) {\Lambda_1},  \nonumber \\
\ k,n \in {\Z}_+ , \ n \le 2u+t-2, \ k \le u-1 \}.  \nonumber  
\eea
Set $$S ^{m} = \{ n -k (m+2) \vert
 k,n \in {\Z}_+ , \ n \le 2u+t-2, \ k \le u-1 \}.$$
\begin{remark}
The weights  ${\l} \in  P^{m}$ was introduced in \cite{KW}.
The corresponding modules are called admissible.
\end{remark}
  
The classification of irreducible $L( m,0)$--modules was
given in \cite{AM}. Recall the following result.
 
\begin{theorem}\label{am} {\bf [AM, Theorem 3.5.3] } 
The set
\bea
\label{oklas}
\{ L(m, r \omega _1) \ \vert \ r \in S ^m \}
\eea
provides all
irreducible $L(m,0)$--modules from the category ${\cal O}$.
\end{theorem}

\begin{remark}
Previous theorem was also proved   in \cite{DLM2}.
\end{remark}

Theorem \ref{am} shows that
 admissible representations of level $m$ for ${\hg}$ can be identified with
the irreducible $L( m,0)$--modules in the category ${\cal O}$.
In the case when $m \in {\N}$ representations defined with (\ref{oklas})
are all irreducible modules for VOA $L(m,0)$. In the case $m \notin {\N}$ and
$m$ admissible, VOA $L(m,0)$ has uncountably many irreducible representations
outside the category ${\cal O}$. The classification of irreducible
$L(m,0)$--modules in the category of weight modules was given in \cite{AM}.
We will now recall the classification result.

For every $r,s \in {\C}$ define
$E_{r,s}=t^{s} {\C}{[t,t^{-1}]}$.
Set $E_{i}=t^{s + i}$. We define $U({\g})$ action on $E_{r,s}$ :
\be \label{form}
e.E_{i}=-(s + i)E_{i-1},\ \ h.E_{i}=(-2 s -2i + r) E_{i}, \ \
f.E_{i}=(s +i-r)E_{i+1} .
\ee
Set $$T^{m}=
\{ \ (r,s)\ \ : \ \ r \in S^{m} \setminus {\Z} , \ \ s \notin {\Z}, \ \
r - s \notin {\Z} \}. $$
Then $E_{r,s}$ is an irreducible $U({\g})$--module  if 
$(r,s) \in T^{m}$.

 Define $\Omega =ef + fe + \frac{1}{2} h^2$ the Casimir element element of
  $U({\g})$.
  The proof of the following lemma is standard.
  \begin{lemma} \label{casimir}
  Let  $w \in E_{r,s}$, $w \in V(r \omega _1)$ or
  $ w \in V(r \omega _1) ^ *$ .Then
  $ \Omega w = \frac {r(r+2)} {2} w $. 
  \end{lemma}
        
We recall the following results from [AM, Section 4].       
\begin{theorem} \cite{AM}    \label{am1}
Assume that $m$ is admissible, $m\notin {\N}$.
Let $r \in S^{m} \setminus {\Z}$.
Then ${\hg}$--module $L(m,E_{r,s} )$ is an  module for the VOA
$L(m,0)$. Moreover, $L(m,E_{r,s} )$ is an irreducible $L(m,0)$--module
if and only if   $(r,s) \in T^{m}$.
\end{theorem}
We need the following consequence of the Theorem \ref{am1}

\begin{corollary}
 \label{tezina}
Assume that $m$ is admissible, $m\notin {\N}$.
Let $r \in S^{m} \setminus {\Z}$.
Then for every $\beta \in {\C}$, there exists an   $L(m,0)$--module $M$
and the weight vector $w \in M_0$ such that
$$ \Omega \vert M_0 = \frac{r(r+2)}{2} \mbox{Id}, \qquad
h(0) w = \beta w. $$
\end{corollary}
{\em Proof.}
Let $s = \frac{ r- \beta} {2}$ and $M= L(m, E_{r,s} )$. Then $M_0 \cong
E_{r,s}$. Set $w = E_0$. Then Lemma \ref{casimir} and relation (\ref{form})
imply  that
$$ \Omega \vert M_0 = \frac{r(r+2)}{2} \mbox{Id}, \qquad
h(0) w = \beta w. {\ \ \Box} $$
                                       
\begin{theorem}\cite{AM}   \label{am.clas}
\label{sl2klas}
 Let $M= \oplus_{n=0}  ^ {\infty }M_n $ be an irreducible $L(m,0)$--module
such that $M_0$ is a weight $U({\g})$--module.
 The $M$ is one of the following
modules:
\bea
&& L(m, V(r \omega _1) ),\quad r \in S^m ,\non \\
&& L(m, V(r \omega _1)^* ),\quad r \in S^m, \non \\
&& L(m,E_{r,s} ), \quad (r,s) \in  T^{m}.   \non
\eea
\end{theorem}

\begin{remark} Modules  $L(m,E_{r,s} )$ are not in the category ${\cal O}$.
In particular, these modules are not irreducible quotient of any Verma modules
over ${\hg}$. In [FST] B. L.  Feigin, A. M. Semikhatov, and I. Yu. Tipunin
introduced the notion of relaxed Verma module.  Our modules
 $L(m,E_{r,s} )$   are   irreducible quotients of relaxed Verma modules.
 \end{remark}

\section{Fermionic and lattice construction of vertex superalgebras}

In this section we will consider vertex superalgebra $F$
constructed from
two charged fermions and vertex superalgebra $F_{-1}$ constructed from
rank one lattice $L= {\Z} \alpha$ such that $ \la \alpha, \alpha \ra = -1$.
In the following sections we will use vertex superalgebras $F$ and
$F_{-1}$ for finding the connections between the representation theory
of VOA $L(m,0)$ and SVOA $L_c$.

 \subsection{Fermionic SVOA $F$}

 Recall the construction of fermionic SVOA $F$.

The charged ree fermionic fields are
\begin{eqnarray*}     
\Psi ^{\pm}(z) = \sum_{i \in \hf + \bf Z} \psi_{i}^{\pm}
z^{-i-\hf},\,\,\,\,\,\, \nonumber
\end{eqnarray*}
with the following  commutation relations
$$ \{ \psi_{i}^{+},\psi_{j}^{-} \}                         
   =  \delta_{i,j},
   \quad
   \{ \psi_{i}^{\pm},\psi_{j}^{\pm} \} =0      
   $$
Let  $F$ be the Fock space defined by
$\psi_{i > 0}^{\pm} {\1} = 0$.  Then $F$ is
a SVOA with the central charge $c=1$ (see \cite{KWn} and \cite{K}
for details).

\subsection{Lattice construction of vertex superalgebras}

  Let $L$  be a  lattice.
  Set ${\h}={\C}\otimes_{\Z}L$ and
extend the ${\Z}$-form $ \la \cdot, \cdot \ra $ on $L$ to ${\h}$.
 Let ${\hh}={\C}[t,t^{-1}]\otimes {\h} \oplus {\C}c$ be the affinization of
${\h}$ (see Section 2).
We also use the notation  $h(n)=t^{n}\otimes h$ for $h\in {\h}, n\in
{\Z}$.

Set 
\begin{eqnarray*}
\hat{{\h}}^{+}=t{\C}[t]\otimes {\h};\;\;{\hh}^{-}=t^{-1}{\C}[t^{-1}]\otimes {\h}.
\end{eqnarray*} 
Then
${\hh}^{+}$ and ${\hh}^{-}$ are abelian subalgebras of
${\hh}$.  
Let $U({\hh}^{-})=S({\hh} ^{-})$ be
the universal enveloping algebra of ${\hh}^{-}$. Consider the
induced ${\hh}$-module
\begin{eqnarray*}
M(1)=U({\hh})\otimes _{U({\C}[t]\otimes {\h}\oplus {\C}c)}{\C}
\simeq S({\hh}^{-})\;\;\mbox{(linearly)},\end{eqnarray*}
where ${\C}[t]\otimes {\h}$ acts trivially on ${\C}$ and $c$ acts
on ${\C}$ as multiplication by 1.

Let
$\hat{L}$ be the canonical central extension of $L$ by the cyclic
group $\< \pm 1\>.$ Form the induced $\hat{L}$-module
\begin{eqnarray*}
{\C}\{L\}={\bf C}[\hat{L}]\otimes _{\< \pm 1\>}{\C}\simeq
{\C}[L]\;\;\mbox{(linearly)},\end{eqnarray*}
where ${\C}[\cdot]$ denotes the group algebra and $-1$ acts on ${\bf
C}$ as multiplication by $-1$. For $a\in \hat{L}$, write $\iota (a)$ for
$a\otimes 1$ in ${\C}\{L\}$. Then the action of $\hat{L}$ on ${\C}
\{L\}$ is given by: $a\cdot \iota (b)=\iota (ab)$ and $(-1)\cdot \iota
(b)=-\iota (b)$ for $a,b\in \hat{L}$.

Furthermore we define an action of ${\h}$ on ${\C}\{L\}$ by:
$h\cdot \iota (a)=\< h,\bar{a}\> \iota (a)$ for $h\in {\h},a\in
\hat{L}$. Define $z^{h}\cdot \iota (a)=z^{\< h,\bar{a}\> }\iota (a)$.

The untwisted space associated with $L$ is defined to be
\begin{eqnarray*}
V_{L}={\C}\{L\}\otimes _{{\C}}M(1)\simeq {\C}[L]\otimes
S({\hh}^{-})\;\;\mbox{(linearly)}.\end{eqnarray*}
Then $\hat{L},{\hh},z^{h}\;(h\in {\h})$ act naturally on
$V_{L}$ by
acting on either ${\C}\{L\}$ or $M(1)$ as indicated above.

For $h\in {\h}$ set
$
h(z)=\sum _{n\in {\Z}}h(n)z^{-n-1}.
$
We use a normal ordering procedure, indicated by open colons, which
signify that in the enclosed expression, all creation operators $h(n)$
$(n<0)$,$a\in \hat{L}$ are to be placed to the left of all
annihilation operators $h(n),z^{h}\;(h\in {\h},n\ge 0)$. For $a \in
\hat{L}$, set
\begin{eqnarray*}
Y(\iota (a),z)=\n e^{\int (\bar{a}(z)-\bar{a}(0)z^{-1})}az^{\bar{a}}\n.
\end{eqnarray*}
Let $a\in \hat{L};\;h_{1},\cdots,h_{k}\in {\h};n_{1},\cdots,n_{k}\in {\Z}\;(n_{i}> 0)$. Set
\begin{eqnarray*}
v= \iota (a)\otimes h_{1}(-n_{1})\cdots h_{k}(-n_{k})\in V_{L}.\end{eqnarray*}
Define vertex operator 
\begin{eqnarray*}
Y(v,z)=\n \left({1\over (n_{1}-1)!}({d\over
dz})^{n_{1}-1}h_{1}(z)\right)\cdots \left({1\over (n_{k}-1)!}({d\over
dz})^{n_{k}-1}h_{k}(z)\right)Y(\iota (a),z)\n 
\end{eqnarray*}
This gives us a well-defined linear map
\begin{eqnarray*}
Y(\cdot,z):& &V_{L}\rightarrow (\mbox{End}V_{L})[[z,z^{-1}]]\nonumber\\ 
& &v\mapsto Y(v,z)=\sum _{n\in {\Z}}v_{n}z^{-n-1},\;(v_{n}\in {\rm
End}V_{L}).
\end{eqnarray*}

Let $\{\;\a_{i}\;|\;i=1,\cdots,d\}$ be an orthonormal basis of ${\bf
h}$ and set
\begin{eqnarray*}
\omega ={1\over 2}\sum _{i=1}^{d}\a_{i}(-1)\a_{i}(-1)\in V_{L}.
\end{eqnarray*}
Then $Y(\omega,z)=\sum_{n\in {\Z}}L_n z^{-n-2}$ gives rise to a
representation of the Virasoro algebra on $V_{L}$ and 
\begin{eqnarray}  \label{vir.rel}
& &L_0\left(\iota(a)\otimes h_{1}(-n_{1})\cdots
h_{n}(-n_{k})\right)\nonumber \\
&=&\left({1\over 2}\< \bar{a},\bar{a}\>+n_{1}+\cdots+n_{k}\right)
\left(\iota(a)\otimes h_{1}(-n_{1})\cdots
h_{k}(-n_{k})\right).
\end{eqnarray}

Now we will assume that $L = {\Z} \alpha$ is a rank one lattice.                                                              
The following theorem is a special case of the
results of Kac [K, Proposition 5.5],
 and Dong and Lepowsky
[DL, Remark  6.17, Remark 9.21].  
\begin{theorem} 
Assume that $L={\Z} \alpha $ is a rank one  lattice,
$\la \alpha, \alpha \ra = n$ and $n$ is add integer.
 Then $V_L$ is vertex superalgebra.
Moreover, if $L$ is a positive definite ( i.e. $n > 0$) then $V_L$ is vertex
operator superalgebra.
\end{theorem}
\begin{remark}
If $n <0$, then  relation
 (\ref{vir.rel}) gives that $V_L$ is a $\frac{1}{2} {\Z}$--graded
 with the respect to $L_0$ and the weight subspaces are not
 bounded below. This implies that
 $V_L$ is  vertex superalgebra which is not vertex operator superalgebra
 (we follow the definitions from \cite{DL} and \cite{Li}).
 \end{remark}
\begin{remark}
If 
$\la \alpha,\alpha \ra=1$, then $V_L$ is SVOA which is isomorphic to the SVOA
$F$ constructed from two charged fermions. This fact is in conformal
field theory known as boson-fermion correspodence.
If  $\la \alpha,\alpha \ra=3$, then $V_L$
is isomorphic to $N=2$ SVOA $L_c$ with $c=1$ (cf. \cite{K}).
\end{remark}

The following discussion is similar as in \cite{DLM}. In \cite{DLM} the authors
considered the case of positive definite even lattice $L$ when $V_L$ is VOA.
We are interested in the  case when $L$ is negative definite lattice
of rank one.

 Define the Schur polynomials $p_{r}(x_{1},x_{2},\cdots)$ $(r\in
{\Z}_{+})$ in variables $x_{1},x_{2},\cdots$ by the following
equation:
\begin{eqnarray}\label{eschurd}
\exp \left(\sum_{n= 1}^{\infty}\frac{x_{n}}{n}y^{n}\right)
=\sum_{r=0}^{\infty}p_{r}(x_1,x_2,\cdots)y^{r}.
\end{eqnarray}
For any monomial $x_{1}^{n_{1}}x_{2}^{n_{2}}\cdots x_{r}^{n_{r}}$ we
have an element $h(-1)^{n_{1}}h(-2)^{n_{2}}\cdots
h(-r)^{n_{r}}{\bf 1}$ in $V_{L}$ for $h\in{\h}.$ Then for any polynomial
$f(x_{1},x_{2}, \cdots)$, $f(h(-1), h(-2),\cdots){\bf 1}$ is a
well-defined element in $V_{L}$.  In particular,
$p_{r}(h(-1),h(-2),\cdots){\bf 1}$ for $r\in {\Z}_{+}$ are elements
of $V_{L}$.

Suppose $a,b\in \hat{L}$ such that $\bar{a}=\alpha,\bar{b}=\beta$. 
Then
\begin{eqnarray}\label{erelation}
Y(\iota(a),z)\iota(b)&=&z^{\<\alpha,\beta\>}\exp\left(\sum_{n=1}^{\infty}
\frac{\alpha(-n)}{n}z^{n}\right)\iota(ab)\nonumber\\
&=&\sum_{r=0}^{\infty}p_{r}(\alpha(-1),\alpha(-2),\cdots)\iota(ab)
z^{r+\<\alpha,\beta\>}.
\end{eqnarray}
Thus 
\begin{eqnarray}\label{eab1}
\iota(a)_{i}\iota(b)=0\;\;\;\mbox{ for }i\ge -\<\alpha,\beta\>.
\end{eqnarray}
Especially, if $\<\alpha,\beta\>\ge 0$, we have
$\iota(a)_{i}\iota(b)=0$ for all $i\in {\Z}_{+}$, and if
$\<\alpha,\beta\>=-n<0$, we get
\begin{eqnarray}\label{eab}
\iota(a)_{i-1}\iota(b)=p_{n-i}(\alpha(-1),\alpha(-2),\cdots)\iota(ab)
\;\;\;\mbox{ for }i\in {\Z}_{+}.
\end{eqnarray}
                              
{\bf From now on we will assume that  $L={\Z} \alpha $ and
$\la \alpha,\alpha \ra=-1$.} 

Set $F_{-1} = V_L$.
Let $a\in \hat L$ such that $\bar a=\alpha.$ 
Set $e=\iota(a), f=\iota(a^{-1}),
k=\alpha(-1){\bf 1}$.  Set $e ^ n = \iota(a^n), f^n = \iota(a^{-n})$.

The relations (\ref{eab1}) and (\ref{eab}) in the case
of vertex superalgebra $F_{-1}$ give the following proposition.
\begin{proposition}   \label{fminus}
 The following relations are hold in the vertex superalgebra $F_{-1}$.
\bea
(a)&& e_i f = 0 \ \ \mbox{for} \ i\ge -1, \quad e_{-2} f = {\1}, \ 
 e_{-3} f = k; \nonumber \\   
(b)&& f_i e = 0 \ \ \mbox{for} \ i\ge -1, \quad f_{-2} e = {\1}, \ 
 f_{-3} e =- k; \nonumber \\   
(c)&& e_i e = 0 \ \ \mbox{for} \ i\ge 1, \quad e_{0} e = e^2, \nonumber \\   
(d)&& f_i f = 0 \ \ \mbox{for} \ i\ge 1, \quad f_{0} f = f^2, \nonumber \\   
(e)&& k_i e = 0 \ \ \mbox{for} \ i\ge 1, \quad k_{0} e = -e, \nonumber \\   
(f)&& k_i f = 0 \ \ \mbox{for} \ i\ge 1, \quad k_{0} f = f. \nonumber    
 \eea
\end{proposition}

\section{ Embedding of $N=2$ SVOA $L_{c}$ into SVOA $F \otimes L(m,0)$}

 The tensor product of the SVOA $F$ and the VOA $L(m,0)$ is SVOA
 $F \otimes L(m,0)$.
 We will show that the SVOA $L_{ c_m }$ can be realized as a subalgebra
 of the SVOA $F \otimes L(m,0)$, where $c_m = \frac{ 3m } {m+2}$.

 Define the following vectors in $F \otimes L(m,0)$:
 \bea
 \tau ^ + =&& \psi_{-{\hf} }^- {\1} \otimes e(-1) {\1} , \label{g+} \\
  \tau ^- =&& \frac{2}{m+2} ,
 \psi_{-{\hf} }^- {\1} \otimes f(-1) {\1}, \label{g-} \\
  j = &&\frac{m}{m+2} \psi_{-\hf}^+ \psi_{- \hf}^- {\1} \otimes {\1}
 - \frac{1}{m+2} {\1} \otimes  h(-1) {\1}, \label{heis} \\
  \nu = &&\frac{1}{m+2} {\1} \otimes e(-1) f(-1) {\1} -
 \frac{m}{m+2}  \psi_{-\hf}^+ \psi_{- \frac{3}{2} }^- {\1} \otimes {\1}
  \non \\
&&- \frac{1}{m+2} \psi_{-\hf}^+ \psi_{- \frac{1}{2} }^-{\1}
 \otimes  h(-1) {\1}. \label{vir}
 \eea
\begin{remark} Our relations (\ref{g+})-(\ref{vir}) are similar to relations
 (3.1)-(3.3) from \cite{FST}.
 \end{remark}

 The following lemma can be proved by direct calculations.
 \begin{lemma}  \label{n2svoa}
 $F\otimes L(m,0)$ is $N=2$ SVOA, i.e. the component of the fields
 $Y(\tau ^+,z)$, $Y(\tau ^-,z)$, $Y(j,z)$ and $Y(\omega,z)$ span $N=2$
 superconformal algebra with the central charge $c = c_m$.
 \end{lemma}
 Now, Corollary \ref{n2con} and Lemma \ref{n2svoa} imply that
 $F \otimes L(m,0)$ is an $U(\cal A)$--module and
 $U(\cal A) . ({\1} \otimes {\1})$ is an subalgebra of the SVOA
 $F\otimes L(m,0)$ isomorphic to certain quotient of the SVOA
 $V_{ c_m}$.

 The mapping
 $U({\cal A}) \rightarrow F \otimes L( m ,0)$ is a special case of
  Kazama-Suzuki
 mapping considered by Feigin, Semikhatov  and Tipunin in \cite{FST}.
 They constructed functor between certain categories
 of ${\hat sl_2}$--modules and $N=2$--modules. They showed that
 this functor is an equivalence of categories.

 Applying this functor in our case we see that
  $U(\cal A) . ({\1} \otimes {\1})$ is an irreducible $U(\cal A)$--module,
  and  it is isomorphic to $L_{ c_m}$. We get the following
   theorem.
       
 \begin{theorem} \label{in2}
 Let $m \ne -2$.
 \item[(1)] SVOA $L_{c_m} \cong U(\cal A) . ({\1} \otimes {\1})$
 is a subalgebra of SVOA $F \otimes L(m,0)$.
 \item[(2)] Assume that $M$ is $L(m,0)$--module. Then $F \otimes M$ is
 a module for the SVOA $L_{ c_m}$.
 \end{theorem}

 \begin{remark}
 The irreducibility of the   submodule
 $U(\cal A) . ({\1} \otimes {\1})$ was also proved in [EG, Appendix A]
    using the calculations of vacuum characters (see also \cite{FSST}).
 \end{remark}

\section{ Embedding of ${ \hat {sl_2} }$ VOA $L(m,0)$
into vertex superalgebra $L_c \otimes F_{-1}$. }
 
 In this section we will consider  tensor product $N=2$ SVOA $L_c$ with
 the lattice vertex superalgebra $F_{-1}$. We will show that the simple
 VOA $L(m,0)$ is a subalgebra of $L_c \otimes F_{-1}$. Our costruction
 is the vertex operator algebra interpretation of the 'anti'-Kazama-Suzuki
 mapping which are considered in [FST].

Let $m \in {\C}$, $m \ne -2$. Recall the definition of
$e, f \in F_{-1}$ from Section 3.

Set
\bea
&& x = G_{-\frac{3}{2}} ^{+} {\1} \otimes f, \label{x} \\
&& y = \frac{m+2}{2} G_{-\frac{3}{2}} ^{-} {\1} \otimes e, \label{y} \\
&& h = - m {\1} \otimes \alpha(-1){\1} + (m+2) T_{-1} {\1} \otimes{\1}.
\label{h}
\eea
\begin{remark}
The relations  (\ref{x})-(\ref{h}) are similar to   relation (3.13) from
 [FST].
 \end{remark}

Set  
$Y(x,z) = \sum _{n \in {\Z} } x(n) z ^{-n-1}$,
$Y(y,z) = \sum _{n \in {\Z} } y(n) z ^{-n-1}$,
$Y(h,z) = \sum _{n \in {\Z} } h(n) z ^{-n-1}$.
Then 
\bea
&& x(n) = \sum _{ i \in {\Z}  } G_{i+ \frac{1}{2} } ^+ \otimes f_{n-i-2},
\label{xn} \\
 && y(n) =\frac{m+2}{2} \sum _{ i \in {\Z}  } 
 G_{i+ \frac{1}{2} } ^- \otimes e_{n-i-2} ,
\label{yn} \\
&& h(n) =-m  {\bf Id} \otimes {\alpha} (n) + (m+2) T_{n} \otimes {\bf Id}.
\label{hn}
\eea

Relations (\ref{xn})-(\ref{hn}) and Proposition \ref{fminus} imply the
following lemma.
\begin{lemma}  \label{sl2n2}
\bea
(a)&& x(n) x = 0, \ \forall n \in {\Z}_+ \qquad y(n) y = 0  \ \forall n\in {\Z}_+,
\nonumber \\
(b)&& x(n) y = 0 \mbox{for} \ n \ge 2, \quad x(1) y = m {\1}, \quad x(0)y=h,
\nonumber    \\
(c)&& h(n) x = 0 \mbox{for} \ n \ge 1, \quad h(0) x = 2 x, 
\nonumber      \\
(d)&& h(n) y = 0 \mbox{for} \ n \ge 1, \quad h(0) y = - 2 y,
\nonumber        \\
(e)&& h(n) h = 0 \mbox{for} \ n \ge 2, \quad h(1)h = 2 m {\1}, \quad h(0)y=0.
\nonumber
\eea
\end{lemma}

Lemma \ref{sl2n2} implies that $L_{ c_m} \otimes  F_{ -1}$ is
a module  for affine Lie algebra ${\hat sl_2}$ at the level $m$.
The mapping
$U({\hat sl_2}) \rightarrow   L_{ c_m} \otimes  F_{ -1}$ is
a spacial case of 'anti'-Kazama-Suzuki mapping considered in
[FST]. This mapping gives a functor from the category
of $N=2$--modules to the category of  ${\hat sl_2}$--modules.
By using  properties of this functor we get that
$U({\hat sl_2}) ({\1} \otimes {\1})$ is an irreducible $U({\hat sl_2})$--module
isomorphic to $L(m,0)$.

We have obtained the following theorem.

\begin{theorem}  \label{ul.sl2}
Let $m \in {\C}$, $m \ne -2$, and
$c= c_m$.
\item[(1)] VOA $L(m,0) \cong U({\hat sl_2}) ({\1} \otimes {\1})$
is a subalgebra of vertex superalgebra \mbox{$L_c \otimes F_{-1}$.}
\item[(2)] Assume that $M$ is a module for SVOA $L_ {c_m}$. Then  $M \otimes F_ {-1}$
 is  a module for VOA $L(m,0)$.

\end{theorem}

\section{Modules for SVOA $L_c$ }

  The representation theory
 of the SVOA $L_c$ with the central charge $c=c_m$ is interesting only in the
 case when $m$ is an admissible rational number.
 Otherwise, $L_c = V_c$ and every highest weight ${\cal A}$--module is
 an $L_c$--module.
 In this section we will present a construction
 of certain set of modules for the SVOA $L_c$ with the central charge
 $c = c_m$. This construction is based on the realization of the SVOA $L_c$
 as a subalgebra of the tensor product SVOA $F \otimes L(m,0)$ and the
 representation theory of the simple VOA $L(m,0)$. When $m \in {\N}$ then
 VOA $L(m,0)$ has finitely many irreducible representations, and we will
 obtain only finitely many irreducible representations for the SVOA
 $L_{ c_m}$. When $ m$ is admissible rational number and $ m\notin {\Q}$
 then VOA $L(m,0)$ has
 uncountably many irreducible representations (see Section 2), and we will get
 uncountably many irreducible $L_c$--representations. 

 For an irreducible $L(m,0)$--module $M$, the top level $M_0$ is an irreducible
 $U(sl_2)$--module.

Now the relations (\ref{g+})-(\ref{vir}) imply the following lemma
(see also [EG, Section 4]).

 \begin{lemma} \label{konstrukcija}
 Let $M$ be any $L(m,0)$--module and $w \in M_0$ such that
 $h(0) w = \beta w$, and
 $\Omega w = \ \gamma w$, for every $w \in M_0$. 
   Then $U({\cal A}) ( {\1}\otimes w)$ is  a highest weight
 $U({\cal A})$--module with the highest weight
 $h,q,c_m$, where
 \be   \label{rel*}
 h= \frac{ \gamma }{2 (m+2)} - \frac{ \beta^2} { 4(m+2)}
  \qquad q= \frac{- \beta}{m+2}.
\ee
 Moreover, $L_{h,q,c_m }$  is an irreducible  $L_{ c_m}$--module.
 \end{lemma}

   For every $r \in {\Z}_+$, $m=\frac{t}{u}$ admissible and $i \in
   \{0,1,\dots,r \}$ we define
  \bea
  &&h_{i,r}= \frac{ r (r+2)}{ 4 (m+2) } - \frac{ (r-2i)^2}{ 4 (m+2) }, \non \\
  && q_{i,r} = \frac{-(r-2i)}{m+2}. \non
  \eea
 Set $N= 2 u + t -2$. Note that if $m \in {\N}$, then $N=m$.
 Define the following finite set
  $$ W^{c_m} = \{
  (h_{i,r},q_{i,r}) \ \vert \   0 \le r \le N, \ 0 \le i \le r \}.
  $$

 Assume now that $M$ is an irreducible $L(m,0)$--module such that
 $M_0$ is finite-\-dimensional. 
  Then $M_0 \cong V(r \omega _1)$ for certain  $r \in \{ 0,\dots,N\}$.
  Then for every weight vector $w \in M_0$
  such that $h(0) w = \beta w$ we have that
   $\beta = r - 2 i$ for certain $i \in
  \{ 0,\dots,r\}$. 
  Now Lemma \ref{konstrukcija} implies the following theorem.
  \begin{theorem} \label{unit.modules}
  Assume that $m \in {\Q}$ is admissible. 
 Then for  every $(h,q) \in W^{c_m} $, the ${\cal A}$--module $L_{h,q,c_m}$
 is the irreducible $L_{c_m}$--module.
 \end{theorem}
\begin{remark} Theorem \ref{unit.modules} shows that starting from
$L(m,0)$--modules $M$ such that $M_0$ is finite-\-dimensional
we can construct only finitely many irreducible modules for SVOA $L_{ c_m}$.
\end{remark}
 We shall now give another parametrisation of the set $W^{c_m}$.

 Let ${\N}_{ \frac{1}{2} }= \{ \frac{1}{2}, \frac{3}{2}, \cdots \}$.
 Set
 $$
  j= i+ \frac{1}{2} \qquad k= r + \frac{1}{2}-i .
 $$
 Then we have that
 $$
 j,k \in {\N} _{\frac{1}{2}  }, \qquad 0 <j ,k, j+k  \le N+1,
 $$
 and   obtain
 $$
 h_{i,r}= \frac{ j k - \frac {1}{4} }{m+2}, \qquad q_{i,r}=\frac {j -k}{m+2}.
 $$
 We get
$$
W^{ c_m} = \left \{ \left( \frac{j k - \frac{1}{4} }{m+2},
 \frac{j-k}{m+2} \right) \ \vert \ j,k \in {\N} _{ \frac{1}{2} },
  0< j,k, j+k \le N+1 \right \},
$$

 In the case $m \in {\N}$ this is exactly the  parametrisation of the
 unitary discrete series
 of $N=2$ minimal models (see \cite{D}).
 So, 
  Theorem \ref{unit.modules}
  shows that every unitary minimal model is a module for SVOA
  $L_{c_m}$.  

Now we assume that $m=\frac{t}{u}$ is admissible and $m \notin {\N}$.
 Then the irreducible
$L(m,0)$--modules are given in    Theorem \ref{sl2klas}.
 Starting from $L(m,0)$--modules $M$ such that
 $M_0$ is finite-dimensional we  constructed finitely many irreducible $L_{ c_m}$--modules
 $L_{h,q,{c_m} }$, $(h,q) \in W^ {c_m}$.
 Now we consider the case when $M_0$ is infinite-\-dimensional.

 \begin{proposition}  \label{nu.modules}  
 Assume that $r \in S^m \setminus {\Z}$, and that $(h,q)$ satisfies the equation
 \be  \label{rel.jed}
 q^2 + \frac{ 4 h}{m+2} = \frac{r (r+2)} { (m+2) ^2}.
 \ee
 Then $L_{h,q,c_m}$ is an irreducible $L_{ c_m}$--module.
 \end{proposition}
 {\em Proof.}
 Assume that $(h,q)$ satisfies the equation (\ref{rel.jed}). Set
 ${\beta} = -q (m+2)$. By using Corollary \ref{tezina} we get that
 there is a $L(m,0)$--module $M$ and the weight vector $ w \in M_0$
 such that
 $$
  \Omega \vert M_0 = \frac{ r (r+2) }{2} \mbox{ Id}, \quad
   h(0) w = \beta w.
$$
Since $(h,q)$ satisfies (\ref{rel.jed}),  we get
$$
q= - \frac{ \beta} { m+2}, \quad
h=  \frac{ \frac{r(r+2)}{2}   } {2 (m+2)} - \frac{ \beta ^ 2} { 4 (m+2)}.
$$
Now Lemma \ref{konstrukcija} implies that $L_{h ,q, c_m}$
 is a $L_{ c_m}$--module. Irreducibility is \mbox {clear. ${\Box}$ }

Define the following set
$$
D^ {c_m} = \left \{ (h,q) \in {\C}^2  \vert \
 q^2 + \frac{ 4 h}{m+2} = \frac{r (r+2)} { (m+2) ^2},
  \ r \in   S^m \setminus {\Z}.      \right \}
$$
Theorem \ref{unit.modules} and Proposition \ref{nu.modules} give
the following theorem.
\begin{theorem}  \label{nonunit.modules}
Assume that $m \in {\Q}$ is admissible and
$(h,q ) \in W^ {c_m} \cup  D^ {c_m}$. Then
$L_{h,q,c_m}$ is an irreducible $L_{c_m}$--module.
\end{theorem}

  \begin{remark}
  Theorem \ref{nonunit.modules} implies that for every admissible nonintegral $m$,
  the SVOA $L_ {c_m}$ has uncountably many non-isomorphic irreducible modules.
  These modules can be parametrized as the union of finitely many rational
  curves in ${\C} ^2$. So, we have constructed uncountably many irreducible modules
  over $N=2$ superconformal algebra which are annihilated with the fields
  $Y(v,z)$, where $v$ is an vector from the maximal submodule of $V_ {c_m}$.
  This gives an difference with the case of Virasoro (cf. \cite{W})
  and Neveu-Schwarz algebra (cf. \cite{A})
  where there exists only finitely many such irreducible modules.
  \end{remark}

\section{Classification of irreducible $L_{c_m}$--modules}

In this section we will give the complete classification of the irreducible
$L_{c_m}$--modules. We will prove that
$L_{ c_m}$--modules constructed in Section 6 gives all irreducible modules
for SVOA $L_{ c_m}$.
The proof of the classification result will use
the fact that VOA $L(m,0)$ is a subalgebra of SVOA
$L_{c_m}\otimes F_{-1}$--modules and the classification of all
irreducible $L(m,0)$--modules obtained in \cite{AM}.

Let ${\g} =sl_2$ and ${\hg} = \hat {sl_2}$ as before.              

\begin{lemma}  \label{clas.l1}
Let $c= c_m$ for $m$ admissible.
Assume that  $L_{h,q,c_m}$ is a $L_{c_m}$--module. Then we have that
$$
\frac{4 h}{m+2} + q^2  = \frac{r(r+2)}{(m+2)^2}
$$
for certain $r \in S^m$.
\end{lemma}
{\em Proof.}
  Assume that $L_{h,q,c}$ is $L_c$--module.
Let $v_{h,q,c}$ be the highest weight vector in $L_{h,q,c}$. Since
$L(m,0)$ is a subalgebra of $L_c \otimes F_{-1}$ (Theorem \ref{ul.sl2}),
we have that
 $L_{h,q,c} \otimes F_{-1}$ is a module for VOA $L(m,0)$--module.
 In particular
  $M= U( {\hg} ) ( v_{h,q,c} \otimes {\1} )$ is an
  $L(m,0)$--module. Set
 $ M_0 = U( {\g}) ( v_{h,q,c} \otimes {\1} )$.
 Since
 $$
 ( {\g} \otimes t^ n )  M_0 = 0 \quad \mbox{for} \ \ n \ge 1,
 $$
  we conclude that the top level of $L(m,0)$--module $M$ is
 $M_0 = U({\g})  ( v_{h,q,c} \otimes {\1} )$.
 Then the representation theory of VOA $L(m,0)$
  (see Theorem \ref{am.clas}) easily implies
  that  $M_0$ is an irreducible  $U({\g})$--module which is
  isomorphic to one of the
  following modules :
  $$
  V(r \omega_1), V(r \omega_1) ^*, E_{r,s} \  \ \mbox{ for} \
  r \in S^m, \ (r,s) \in T^m.
  $$
 Let $\Omega = x(0)y(0) + y(0)x(0) + \frac{1}{2}h(0)^2$
be the Casimir.
Then Lemma  \ref{casimir}
imply that
\be \label{cas2}
\Omega \vert M_0 = \frac{r(r+2)}{2}  \mbox{ Id} \quad \mbox{for certain} \ r\in S ^m.
 \ee

Let $w \in M_0$. Then there is $f \in
U({\g})$ such that $ w = f ( v_{h,q,c} \otimes {\1} )$.
 Since the action of $\Omega$
 commutes with $U({\g})$, 
 then from relations (\ref{xn})-(\ref{hn}) we get
\bea
&& \Omega w =f ( \Omega ( v_{h,q,c} \otimes {\1} ) ) \nonumber \\
&&= f (x(0)y(0) + y(0)x(0) + \frac{1}{2}h(0)^2) ( v_{h,q,c} \otimes {\1} )
\nonumber \\
&& =  ( 2(m+2)h + \frac{1}{2}(m+2)^2 q^2 )   w.
\eea
So, we have proved 
\be \label{cas1} 
\Omega \vert M_0 = ( 2(m+2)h + \frac{1}{2}(m+2)^2 q^2 ) \mbox{ Id}.
\ee
 Now from (\ref{cas2}) and (\ref{cas1}) follow that
$$
\frac{4 h}{m+2} + q^2  = \frac{r(r+2)}{(m+2)^2}
$$
for certain $r \in S^m$.  ${\ \ \Box}$

\begin{lemma} \label{clas.l2}
 Assume that $L_{h,q,c_m}$ is $L_{c_m}$--module and
\be     \label{rl}
\frac{4 h}{m+2} + q^2  = \frac{r(r+2)}{(m+2)^2}
\ee
for $r \in S^m \cap {\Z}_+$. Then $(h,q) \in W^{ c_m}$.
\end{lemma}
{\em Proof.}  Let $M= U( {\hg} ) ( v_{h,q,c} \otimes {\1})$ be a
$L(m,0)$--module as
in the proof of Lemma \ref{clas.l1}.
Since in the relation (\ref{rl}) we have $r \in {\Z}_+$, we conlclude that
$M_0$ is an irreducible finite-dimensional $U({\g})$--module isomorphic to
$V(r \omega _1)$. This implies that the vector
$ v_{h,q,c} \otimes {\1} $ is a weight vector of $V(r \omega _1)$,
i.e.
$$
h(0) . (v_{h,q,c} \otimes {\1}) = (r-2i)(v_{h,q,c} \otimes {\1})
$$
 for certain
$i \in \{0,\dots,r\}$. Applaying the formulae (\ref{hn}) for the action
of $h(0)$ on $L_{h,q,c} \otimes F_{-1}$  
 we get
$   q=\frac{r-2i}{m+2}$. Now the relation (\ref{rl}) implies that
 $h=   \frac{r(r+2)}{4 (m+2) } - \frac{ (r-2i)^2 }{4 (m+2) }$, and
 we conclude that $(h,q) \in W^{ c_m}$.
 ${\ \ \Box}$

\begin{theorem} \label{unit.classification}
Assume that $m \in {\N}$ and $c = c_m$.  Then the set
$$
\{L_{h,q,c} \ \vert \  (h,q) \in W^ {c_m} \}
$$ provides all $L_c$ irreducible modules for
the SVOA $L_c$.
So, irreducible $L_c$--modules are exactly all unitary modules for
$N=2$ superconformal algebra with the central charge $c$.
\end{theorem}
{\em Proof.}
We proved in Theorem \ref{unit.modules} that for every  $(h,q) \in W^ {c_m}$,
$L_{h,q,c}$ is a $L_c$--module.  It remains to prove that if
$L_{h,q,c}$ is a $L_c$--module, then  $(h,q) \in W^ {c_m}$.

Assume now that  $L_{h,q,c}$ is a
 $L_c$--module.
  Then Lemma \ref{clas.l1} implies that
$$ 
\frac{4 h}{m+2} + q^2  = \frac{r(r+2)}{(m+2)^2}
$$ 
for certain $r \in S^m$.
Since $S^m \subset {\Z}_+$ for $m \in {\N}$, we
have that $r \in {\Z}_+$. Now Lemma \ref{clas.l2} implies that
$(h,q) \in  W^ {c_m}$. ${\ \ \Box}$

\begin{remark}
Theorem \ref{unit.classification} shows that SVOA $L_{c_m}$ for $m \in {\N}$
has exactly $\frac{(m+2) (m+1)} {2} $ non-isomorphic irreducible modules.
\end{remark}
\begin{theorem} \label{clas.nonunit}
Assume that $m \in {\Q}$ is admissible such that $ m \notin {\N}$. Let 
$c = c_m$. Then the
set
$$
\{ L_{h,q,c} \ \vert \  (h,q) \in  W^ {c_m} \cup D^ {c_m} \}
$$
provides all $L_c$ irreducible modules for the 
SVOA $L_c$.
\end{theorem}
{\em Proof.}
Theorem \ref{nonunit.modules} gives  that $L_{h,q,c}$ is a $L_c$--module
for every
\mbox{ $ (h,q) \in  W^ {c_m} \cup D^ {c_m}$}. In order
to prove theorem we have to prove that if $L_{h,q,c}$ is a $L_c$--module,
then  $ (h,q) \in  W^ {c_m} \cup D^ {c_m}$.

Assume now that $L_{h,q,c}$ is a $L_c$--module.
  Then Lemma \ref{clas.l1} implies that
$$ 
\frac{4 h}{m+2} + q^2  = \frac{r(r+2)}{(m+2)^2}
$$
for certain $r \in S^m$.
If $r \in   S^m \setminus {\Z}$, then $(h,q) \in D^ {c_m}$.
If $r \in {\Z}_+$, then from Lemma \ref{clas.l2} follows that
$(h,q) \in W^ {c_m}$. So, we get  
\mbox{$ (h,q) \in  W^ {c_m} \cup D^ {c_m}$.   ${ \ \Box}$ }

Department of Mathematics, University of Zagreb, 
Bijeni\v{c}ka 30, 10000 Zagreb, Croatia

E-mail address: adamovic@cromath.math.hr

\end{document}